\documentclass[12pt]{article}
\usepackage{amssymb,amsmath, amsthm, amscd,ifthen}
\usepackage[T2A]{fontenc}
\usepackage[cp1251]{inputenc}
\usepackage[english,russian]{babel}
\usepackage[dvips]{graphicx}
\usepackage[hyper]{amsbib}
\usepackage{indentfirst}

\newtheorem{theorem}{Theorem}

\newtheorem{proposition}{Proposition}
\newtheorem{lemma}{Lemma}
\newtheorem{corollary}{Corollary}

\newenvironment{proof1}{\noindent{\it Proof.\,}}{\hfill$\Box$}

\begin{document}

\def\q#1.{{\bf #1.}}
\def\N{\mathbb N}
\def\Z{\mathbb Z}
\def\Q{\mathbb Q}
\def\R{\mathbb R}
\def\Pp{\mathbb P}
\def\ffi{\varphi}
\def\vec{\overrightarrow}

\title{Coloring hypergraphs with bounded cardinalities of edge intersections}

\author{Margarita Akhmejanova\footnote{Moscow Institute of Physics and Technology, Laboratory of Advanced Combinatorics and Network Applications, 141700, Institutskiy per. 9, Dolgoprudny, Moscow Region, Russia. E-mail: mechmathrita@abc.math.msu.su}, Dmitry Shabanov\footnote{Moscow Institute of Physics and Technology, Laboratory of Advanced Combinatorics and Network Applications, 141700, Institutskiy per. 9, Dolgoprudny, Moscow Region, Russia; National Research University Higher School of Economics (HSE), Faculty of Computer Science, 101000, Myasnitskaya Str. 20, Moscow, Russia. E-mail: dmitry.shabanov@phystech.edu.}}
\date{}
\maketitle

\bigskip
\textbf{Abstract.} The paper deals with an extremal problem concerning colorings of hypergraphs with bounded edge degrees. Consider the family of $b$-simple hypergraphs, in which any two edges do not share more than $b$ common vertices. We prove that for $n\geqslant n_0(b)$, any $n$-uniform $b$-simple hypergraph with the maximum edge degree at most $c\cdot nr^{n-b}$ is $r$-colorable, where $c>0$ is an absolute constant. We also establish some applications of the main result.

\section{Introduction}\label{intro}
The paper deals with colorings of uniform hypergraphs. Let us start with recalling some definitions.

\subsection{Definitions}
A \emph{vertex $r$-coloring} of a hypergraph $H=(V,E)$ is a mapping from the vertex set $V$ to the set of $r$ colors $\{1,\ldots,r\}$. A coloring of $H$ is called \emph{proper} if there is no monochromatic edges under this coloring, i.e. every edge of $H$ contains at least two vertices which receive different colors. A hypergraph is said to be $r$-\emph{colorable} if there exists a proper $r$-coloring for it. The \emph{chromatic number} of hypergraph $H$ is the minimum $r$ such that $H$ is $r$-colorable.

The \emph{degree of an edge} $A$ in a hypergraph $H$ is the number of other edges of $H$ which have nonempty intersection with $A$. The maximum edge degree of $H$ is denoted by $\Delta(H)$. For a given natural number $b$, a hypergraph $H=(V,E)$ is said to be $b$-\emph{simple} if every two distinct edges of $H$ do not share more than $b$ common vertices, i.e., formally,
$$
  |A\cap B|\leqslant b\mbox{ for any }A,B\in E,\;\;A\ne B.
$$
The main aim of the current work is to refine a quantitative relation between the chromatic number and the maximum edge degree in an $n$-uniform $b$-simple hypergraph.

\subsection{Related work}
The first quantitative relation between the chromatic number and the maximum edge degree in a uniform hypergraph was obtained by Erd\H{o}s and Lov\'asz in their classical paper \cite{ErdLov}.
They proved that if $H$ is an $n$-uniform hypergraph and
\begin{equation}\label{bound:erdlov}
  \Delta(H)\leqslant \frac 14r^{n-1},
\end{equation}
then $H$ is $r$-colorable. Recall that the result was historically the first application of the Local Lemma. The bound \eqref{bound:erdlov} appeared not to be sharp. The restriction on the maximum edge degree which guarantees $r$-colorability have been successively improved in a series of papers. We will mention only the best known results, the reader is referred to the survey \cite{RaigShab} for the detailed history of the question.

In connection with the ``Property B'' problem Radhakrishnan and Srinivasan \cite{RadhSrin} proved that any $n$-uniform hypergraph $H$ with
$$
  \Delta(H)\leqslant 0,17 \left(\frac n{\ln n}\right)^{\frac 12}2^{n-1}
$$
is $2$-colorable. The complete generalization of the above result to an arbitrary number of colors was derived by Cherkashin and Kozik \cite{CherKozik}, who showed that any $n$-uniform hypergraph $H$ satisfying
\begin{equation}\label{bound:cherkozik}
  \Delta(H)\leqslant c(r)\left(\frac{n}{\ln n}\right)^{\frac {r-1}{r}}r^{n-1}
\end{equation}
is $r$-colorable, where $c(r)>0$ does not depend on $n$ and $n>n_0(r)$ is large enough. The results \eqref{bound:cherkozik} was derived by the help of Pluh\'ar's approach to colorings of hypergraphs from \cite{Pluhar}. Recent advances concerning ``Property B''-type problems on colorings of hypergraphs can be found in \cite{KupShab}, \cite{Cherk}, \cite{AkhShab}.

\bigskip
In the current work we concentrate on the similar problem in the class of $b$-simple hypergraphs. The case of 1-simple hypergraphs (known also as \emph{linear} or \emph{simple} hypergraphs) is very well studied. The reader is referred to the survey \cite{RaigShab} and the papers \cite{KozikShab}, \cite{KupShab} for the detailed history of the question, we will note only the best current results. Kozik and Shabanov \cite{KozikShab} proved that any $n$-uniform simple hypergraph $H$ with
\begin{equation}\label{bound:kozikshab}
  \Delta(H)\leqslant c\cdot n\, r^{n-1}
\end{equation}
is $r$-colorable, here $c>0$ is some absolute constant. In the case of very large number of colors a better estimate was derived by Frieze and Mubayi \cite{FriezeMubayi}. They showed that the inequality
$$
    \Delta(H)\leqslant c(n)r^{n-1}\ln r
$$
implies $r$-colorability of an $n$-uniform simple hypergraph $H$, where $c(n)>0$ is a very small function of $n$ (the calculations in \cite{FriezeMubayi} give $c(n)\leqslant n^{-2n}$).

\bigskip
The first result concerning colorings of $b$-simple hypergraph in general situation was obtained by Kostochka and Kumbhat \cite{KostKumb}. They proved that for any $\varepsilon>0$, $b\geqslant 1$ and $r\geqslant 2$, arbitrary $n$-uniform $b$-simple hypergraph $H$ with
\begin{equation}\label{bound:kk}
  \Delta(H)\leqslant n^{1-\varepsilon}r^{n-1}
\end{equation}
is $r$-colorable provided $n>n_0=n_0(r,b,\varepsilon)$ is large enough. Since $\varepsilon>0$ is arbitrary in \eqref{bound:kk} then, of course, it can be replaced by some infinitesimal function $\varepsilon=\varepsilon(n)>0$, tending to 0 with growth of $n$. Several papers were devoted to estimating the order of $\varepsilon(n)$. Kostochka and Kumbhat stated that one can take $\varepsilon(n)=\Theta(\frac{\ln\ln\ln n}{\ln\ln n})$. Shabanov \cite{Shab} refined this to $\varepsilon(n)=\Theta((\frac{\ln\ln n}{\ln n})^{1/2})$. The best known result was obtained by Kozik \cite{Kozik}, who showed that for sufficiently large $n$, any $b$-simple $n$-uniform hypergraph $H$ with
\begin{equation}\label{bound:kozik}
  \Delta(H)\leqslant c\cdot\frac{n}{\ln n}r^{n-b-1}
\end{equation}
is $r$-colorable, wherein $c>0$ is some absolute constant.

\subsection{Main result}

The main result of the paper improves the estimate \eqref{bound:kozik} as follows.

\begin{theorem}\label{thm:main}Suppose $b\geqslant 1$, $r\geqslant 2$ and $n>n_0(b)$ is large enough in comparison with $b$. Then if a $b$-simple $n$-uniform hypergraph $H=(V,E)$ satisfies the inequality
\begin{equation}\label{bound:new}
  \Delta(H)\leqslant \frac 1{(2e)^{4}}\cdot n\,r^{n-b},
\end{equation}
 then $H$ is $r$-colorable.
\end{theorem}

In the case of simple hypergraphs, for $b=1$, the result \eqref{bound:new} coincides with \eqref{bound:kozikshab}. Note that for fixed $r,b$, the bound \eqref{bound:new} is at most $n$ times smaller than the best possible. Recall that Kostochka and R\"odl \cite{KR} showed that there exists an $n$-uniform non-$r$-colorable simple hypergraph $H$ with $\Delta(H)\leqslant n^2r^{n-1}\ln r$.

\bigskip
The remaining paper is structured as follows. Section 2 is devoted to the proof of Theorem \ref{thm:main}. In Section 3 we deduce few corollaries.

\section{Proof of Theorem \ref{thm:main}}
The proof of Theorem \ref{thm:main} is based on the random recoloring method. We use its modification from the paper of Kozik and Shabanov \cite{KozikShab}. However we had to derive some new ideas and constructions for application of this method to the case of $b$-simple hypergraphs.

\subsection{The algorithm of recoloring method}

The general principle of the recoloring approach is clear: for a given non-proper coloring of a hypergraph vertex set, we try to recolor a small number to vertices to make the coloring proper.

Suppose $H=(V,E)$ is a $b$-simple $n$-uniform hypergraph satisfying the condition \eqref{bound:new}. We are going to use the randomized algorithm from \cite{KozikShab} to find a proper coloring with $r$ colors for $H$. Let us describe it.

\begin{enumerate}
  \item Consider a random $r$-coloring $f=(f(v),v\in V)$ of the vertex set with uniform distribution on $\{0,\ldots,r-1\}^{|V|}$.
  \item For every vertex $v\in V$, consider an independent random variable $\sigma(v)$ with uniform distribution on $[0,1]$ (also independent of $f$). The value $\sigma(v)$ is called \emph{the weight} of the vertex $v$. With probability $1$ the mapping $\sigma: V\to [0,1]$ is injective.
  \item Given parameter $p$, a vertex $v$ is said to be \emph{free} if $\sigma(v)\leqslant p$. Only free vertices are allowed to be recolored during the recoloring process.
  \item Starting with $f$, do the following.

  \emph{\textbf{Recoloring step}}. If there exists a monochromatic edge $A$ in the current coloring whose first (i.e. a vertex with the least weight) non-recolored vertex $v$ is free then recolor $v$ with color $(f(v)+1) ({\rm mod}\; r)$. In the above situation we say that a vertex $v$ \emph{blames} an edge $A$.
  \item Repeat the recoloring step until possible.
\end{enumerate}

Note that every vertex can be recolored only once during the recoloring procedure, so the process always stops.

\bigskip
Let us understand what configuration can be guilty of failure of the algorithm. Kozik and Shabanov showed the following tree-type construction should take place.

\subsection{H-tree construction}

Suppose that the algorithm fails to produce a proper coloring and an edge $A$ is monochromatic in the final coloring. In the process of considering how the edge $A$ became monochromatic we will build a labelled rooted tree graph, for which we will use the term ``h-tree'' in our text.

\bigskip
Before we begin, let us remind basic facts which follow from the recoloring algorithm. If during the evaluation of the algorithm some vertex $v$ is recolored, then it should be the first non-recolored free vertex of some edge $F$ that at that moment of the procedure is monochromatic, i.e. the vertex $v$ blames the edge $F$. In the case of multiple set of blamed edges we choose one for every vertex. Note that every edge can be blamed only by one vertex.

\bigskip
The construction of h-tree can be build as follows:
\begin{itemize}

 \item An edge $A$ is monochromatic in the final coloring.
  \emph{Create a graph $T$ consisting of a single root-node labelled by the edge $A$. }

  \item Let $a$ be the color of $A$. Edge $A$ is monochromatic in the final coloring so it cannot have free vertices with initial color $a$. Therefore it can contain only non-free vertices with initial color $a$ and free vertices with initial color $a-1$. Vertices of the last type blame some edges, say, $B_1,\ldots,B_t$. In this situation we also say that $A$ \emph{blames} the edges $B_1,\ldots,B_t$. \emph{Add $t$ new nodes labelled by $B_1,\ldots,B_t$ as children of root-node labelled by $A$.}

  \item Since every $B_i$ should be monochromatic of color $a-1$ at some step of the recoloring procedure then in the initial color it contains only vertices of colors $a-1$ and $a-2$. All the vertices of initial color $a-2$ should be free and should blame some other edges $C^i_1,\ldots,C^i_{t_i}$.\emph{ Add $t_i$ new nodes labelled by  $C^i_1,\ldots,C^i_{t_i}$ as children of node labelled by $B_i$.}
  \item Continue the process until possible.
\end{itemize}

The obtained configuration $T$ has a tree construction, its vertices, called \emph{nodes} are labelled by edges of $H$, futhermore, its leaves are labelled by the edges that are monochromatic in the initial coloring $f$. (To reduce notation we will denote labelling function by one symbol $\phi:V(T)\rightarrow E(H)$ and will say that node $u$ is labelled by the edge $\phi(u)$ or $\phi(u)$ is the label of $u$.) Also note that the adjacency in $T$ is induced by the blaming relationship ($B$ is a child of $A$ in $T$ if and only if the edge $\phi(A)$ blames edge $\phi(B)$ in $H$). However we will draw the attention to the fact that two different nodes in $T$ can be labelled by the same edge.

Proposition 4 in \cite{KozikShab} states that if the recoloring algorithm fails then there exists an h-tree.

\subsection{The Local Lemma}
The authors of \cite{KozikShab} used a specific variant of the Local Lemma. It is derived from the general version by Beck in \cite{Beck}.

\begin{lemma}\label{LocalLemma}{\rm (Local Lemma)}
Let $\mathcal{X}=\{X_{1},X_{2}\ldots,X_{m}\}$  be independent random variables (or vectors) in arbitrary probability space and let $\mathcal A$ be a finite set of events determined by these variables. For $A\in\mathcal A$, let {\rm vbl}(A) denote the set of variables that determines $A$ (i.e. $A$ belongs to the sigma-algebra generated by $Y\in {\rm vbl}(A)$). For $X\in\mathcal{X}$, define a polynomial $w_{X}(z)$ as follows:
\begin{equation}\label{local:polynomial}
w_{X}(z)=\sum_{A\in\mathcal A:X\in {\rm vbl}(A)}{\sf Pr}(A)z^{|{\rm vbl}(A)|}.
\end{equation}
Suppose that a polynomial  $w(z)$ dominates all the polynomials $w_{X}(z)$ i.e. for every real $z_{0}\ge 1$ we have  $w(z_{0})\geqslant w_{x}(z_{0})$. If there exists $\tau_{0}\in (0,1)$ such that for every $X\in\mathcal{X}$,
\begin{equation}\label{local:condition}
w_X\left(\frac{1}{1-\tau_{0}}\right)\leqslant \tau_{0},
\end{equation}
then all the events from $\mathcal A$ can be simultaneously avoided with positive probability, i.e.
$$
    {\sf Pr}\left(\bigcap_{A\in\mathcal A}\overline{A}\right)>0.
$$
\end{lemma}

The proof of Lemma \ref{LocalLemma} can be found in \cite{Kozik}. In our model the independent random vectors $(f(v),\sigma(v))$, $v\in V$, are labelled by the vertices of the hypergraph. We will estimate the probabilities of the bad events and then sum up those of them (with coefficients in \eqref{local:polynomial}) for which the corresponding bad configuration contains an arbitrary fixed vertex $v$. The choice of the parameters will be the following:
\begin{equation}\label{choice:parameters}
    \tau_0=\frac 1{n+1},\;\;p=\frac{5\ln n}{n}.
\end{equation}
Now we proceed to the analysis of the bad events.

\subsection{Analysis of the bad events}

Suppose that the randomized recoloring algorithm fails. Let $A$ denote the monochromatic edge in the final coloring and let $T$ denote an $h$-tree with root is labelled by $A$.

\subsubsection{Bad event 1: a lot of recolored vertices}
The first bad event $\mathcal{B}_1$ happens if there is an edge $F$ with at least $20e\ln n$ vertices recolored during the procedure. This event implies that
\begin{itemize}
  \item $F$ becomes monochromatic of some color $\alpha$ during the recoloring procedure;
  \item every vertex $v\in F$ either has initial color $f(v)=\alpha$ or $f(v)=\alpha-1\;({\rm mod}\;r)$;
  \item the number of free vertices in $F$ (number $k$) is at least $20e\ln n$;
  \item all the vertices with the initial color $\alpha-1$ are free.
\end{itemize}
Let $\mathcal{B}_1(F)$ denote the event described above. Its probability can be easily calculated:
\begin{align*}
    {\sf Pr}(\mathcal{B}_1(F))&=r\sum_{k\geqslant 20e\ln n}{n\choose k}\left(\frac {1}r\right)^{n-k}\left(\frac{2p}r\right)^{k}=
    r^{1-n}\sum_{k\geqslant 20e\ln n}{n\choose k}(2p)^k\leqslant\\
    &\leqslant r^{1-n}\sum_{k\geqslant 20e\ln n}\left(\frac{2enp}k\right)^k=r^{1-n}\sum_{k\geqslant 20e\ln n}\left(\frac{10e \ln{n}}k\right)^k\leqslant\\
    &\leqslant r^{1-n}\sum_{k\geqslant 20e\ln n}\left(\frac 12\right)^k\leqslant r^{1-n}2^{1-20e\ln n}\leqslant 2\,r^{1-n}n^{-10}.
\end{align*}
Therefore, for every vertex $v$, the following estimate for the local polynomial holds:
\begin{align}\label{polynom:event1}
    w^1_v\left(\frac 1{1-\tau_0}\right)&=\sum_{F:\;v\in F}{\sf Pr}(\mathcal{B}_1(F))\left(\frac 1{1-\tau_0}\right)^{|F|}=\notag\\
&\mbox{(since the number of edges containing $v$ does not exceed $\Delta(H)+1\leqslant 2\Delta(H)$)}\notag\\
    &=\sum_{F:\;v\in F}{\sf Pr}(\mathcal{B}_1(F))\left(1+\frac 1n\right)^{n}\leqslant 2\Delta(H)\cdot 2\,r^{1-n}n^{-10} e\leqslant\notag\\
&\mbox{(using condition \eqref{bound:new})}\notag\\
    &\leqslant \frac{2}{(2e)^4}\,r^{n-b}n\cdot 2\,r^{1-n}n^{-10} e=\frac 1{4e^3} r^{1-b}n^{-9}\leqslant \frac{1}{10(n+1)}.
\end{align}
A label $F$ is called \emph{degenerate} if the event $\mathcal{B}_1(F)$ holds. Now we will consider $h$-trees without degenerate labels.

\subsubsection{Removing the coinciding edges}

Suppose now $T$ is an h-tree with root labelled by $A$ and without degenerate labels. For any node $C\in T$, let $N(C)$ denote a set of all descendants of $C$ in $T$, i.e. $N(C)$ consists of all the nodes $B$ such that $C$ lies on the shortest path from $B$ to $A$ in $T$. The induced subgraph on $N(C)$ forms an \emph{h-subtree}.

The adjacency in  h-tree is induced by the blaming relationship. If node $C$ is a child of node $B$ then $\phi(B)$ contains a vertex $v(\phi(C))$ which blames $\phi(C)$, so $v(\phi(C))\in \phi(C)\cap \phi(B)$.

\bigskip
The next proposition says that every blaming vertex is uniquely defined.
\begin{proposition}\label{blaming}
Suppose that nodes $F_1,\ldots,F_s$ are  children of node $C$ in $h$-tree. Then there is a one-to-one correspondence between the set of edges $\phi(F_1),\ldots,\phi(F_s)$ and the set of blaming vertices $v(\phi(F_1)),\ldots,v(\phi(F_s))$.
\end{proposition}
\begin{proof1}
Edge $\phi(C)$ becomes monochromatic of a color $\alpha$ at some moment of the recoloring procedure. Therefore, there is the last recolored vertex $u$. The vertex $u$ was recolored since it blamed some $\phi(F_i)$, and on that recoloring step $\phi(F_i)$ was completely monochromatic of color $\alpha-1$. Hence, $|\phi(C)\cap \phi(F_i)|=1$ and $u=\phi(C)\cap \phi(F_i)=v(\phi(F_i))$. Let us remove the vertex $u$ from $\phi(C)$ and repeat the above argument. We will obtain the complete one-to-one correspondence between $\phi(F_1),\ldots,\phi(F_s)$ and $v(\phi(F_1)),\ldots,v(\phi(F_s))$.
\end{proof1}

\bigskip
Suppose that $C$ and $D$ are different nodes of $T$, but for their labels it holds that $\phi(C)=\phi(D)$. This situation appears when the blaming vertices, $v(\phi(C))$ and $v(\phi(D))$, coincide. Note that if the labels of the nodes $C$ and $D$ coincide then the same holds for the corresponding nodes in the subtrees $N(C)$ and $N(D)$, so this property is hereditary. So, we say that a vertex $v$ is \emph{special} if there are two different nodes $C$ and $D$ in $T$ such that
\begin{itemize}
  \item $\phi(C)=\phi(D)$;
  \item $v=v(\phi(C))$ and $v=v(\phi(D))$;
  \item Parents of $C$ and $D$ have different labels.
\end{itemize}
The notation of a special vertex in an h-subtree is defined in the same way.

\bigskip
For given h-tree (or h-subtree) $T$, let us define an operation of removing nodes with coinciding labels.
\begin{enumerate}
  \item Let us fix an arbitrary order $\zeta'$ of the edges of $H$. Let us order the nodes of the h-tree $T$ as follows: order them in increasing distance from the root, if the distance is the same then order them according to $\zeta'$, if both, distance and $\zeta'$, are the same (i.e. we have the pair of coinciding labels) then order according to the order of the parents in the h-tree. Let $\zeta$ denote the obtained order.
  \item Consider the nodes according to $\zeta$.
  \item For the current node $C$, if there is a  node $D$ labelled by $\phi(C)=\phi(D)$, copy of $C$, then remove from the h-tree all the copies of $C$ together with all the descendants (i.e. remove $N(D)$ if $D$ is a copy).
  \item Repeat the previous step until possible.
\end{enumerate}
Let $O(T)$ denote the obtained h-tree (or h-subtree). Now $O(T)$ does not contain any nodes with coinciding labels, so we will call it \emph{proper}.

\subsubsection{Bad event 2: $b$-disjoint proper h-trees}

Suppose now $T_1=O(T)$ is a proper h-tree with root labelled by $A$ and without degenerate labels. The key notation of our probabilistic analysis is the notation of a bad node in an h-tree. A node $C\in T$ of the h-tree $T_1$ is said to be \emph{bad} if
$$
  \left|\phi(C)\cap \bigcup_{B\in T_1\setminus N(C)}\phi(B)\right|\geqslant b+1,
$$
i.e. $\phi(C)$ as an edge of the initial hypergraph $H$ has a lot of common vertices with the union of all other edges that are not labels of its descendants in $T_1$. The definition of a bad node in a proper h-subtree is absolutely the same.

\bigskip
\emph{A straight path} in a rooted tree is the shortest path that joins a node of the tree with the root. Now we have two possibilities: either there is a straight path that contains all the bad nodes in the configuration or there is no such path. Suppose the first alternative holds and let $C_m$ denote the bad node with the largest distance from the root $A$. Let $(C_m,C_{m-1},\ldots,C_0=A)$ denote the straight path from $C_m$ to $A$, so all the bad nodes are contained in this set. Assume also that for every $j=0,\ldots,m-1$,
\begin{equation}\label{condition:path}
    \left|\phi(C_j)\cap \bigcup_{i=j+1}^{m}\phi(C_i)\right|\leqslant b.
\end{equation}
A proper h-tree (or h-subtree) is said to be $b$-\emph{disjoint} if there is unique straight path containing all the bad nodes and the condition \eqref{condition:path} holds. Note that if there is no bad nodes then we assume that the path consists only of the root.
In the second bad event $\mathcal{B}_2(T_1)$ we consider the case of $b$-disjoint proper h-tree $T_1$. Let $t$ be the size of $T_1$. The following statements describe the properties of $T_1$.

\begin{proposition}\label{bdisjoint:property1}
The total number of hypergraph vertices over the  all labels of $T_1$ (i.e. in $\phi(T_1)$) is at least $n+(n-b)(t-1)$.
\end{proposition}
\begin{proof1} Let $(C_m,\ldots,C_0)$ denote the path containing all the bad nodes. Let us arrange all the nodes of $T_1$ in the following way: $C_m,\ldots,C_0$ and then all the remaining nodes are ordered according to the distance from the root $C_0$ in the increasing order. Then every node (except $C_m,\ldots,C_0$) is ordered before all its descendants, so its label has at most $b$ common vertices with all labels of the preceding nodes. Otherwise this node is bad, which contradicts the $b$-disjoint property. Due to \eqref{condition:path} the same is true for $C_0,\ldots,C_m$. The label $\phi(C_m)$ consists of $n$ vertices and all labels of the other nodes give at least $n-b$ new vertices, not contained in  the preceding label.
\end{proof1}

\bigskip
Recall that $T_1$ does not have degenerate labels. So, in every label $\phi(C)$ there are at most $20e\ln n$ vertices which were recolored before $\phi(C)$ became monochromatic of some color $\alpha$ in the recoloring procedure. Let us call $\alpha$ the \emph{dominating} color of $C$. At least $n-20e\ln n$ are colored with $\alpha$ in $\phi(C)$ in the initial coloring $f$.

The h-tree construction guarantees that for a fixed dominating color $\alpha$ for $A$, all the colors of the vertices in $\phi(T)$ are uniquely defined in the initial coloring $f$. This immediately follows from the algorithm. Proposition \ref{bdisjoint:property2} says that the same holds for proper h-trees also.

\begin{proposition}\label{bdisjoint:property2}
For a given dominating color of $A$, the initial colors of all the vertices in $\phi(T_1)=\phi(O(T))$ are uniquely defined.
\end{proposition}
\begin{proof1}
For a given dominating color $A$, the dominating colors of all the nodes are uniquely defined. The removing procedure (operation $O$) implies that all the children of the root $A$ cannot be removed. So, due to Proposition \ref{blaming} the set of blaming vertices $v(\phi(F_1)),\ldots,v(\phi(F_s))$ can be uniquely defined by the $\phi(F_1),\ldots,\phi(F_s)$. So, the colors in $A$ are restored.

Now, order the remaining nodes according to $\zeta$ and set $\mathcal{R}(T_1)=\{v(\phi(F_1)),\ldots,v(\phi(F_s))\}$. Then for every next node $C$,
\begin{itemize}
  \item we know all the colors of the common vertices $\phi(C)$ with labels of its children, add of them to $\mathcal{R}(T_1)$;
  \item we know all the colors in the set $\mathcal{R}(T_1)\cap \phi(C)$;
  \item the initial color of all the remaining vertices is the dominating color of $C$.
\end{itemize}
Indeed, if the initial color of a vertex $w$ is not equal to the dominating color of $C$ then this vertex blames another edge $\phi(D)$, where in $D$ is a children of $C$ in the h-tree. Maybe $D$ was removed by the operation $O$, but in this case $T_1$ contains a node $D'$, a copy of $D$. The parent of $D'$ has the less order in $\zeta$ than $C$ and its label contains $w$, so $w\in\mathcal{R}(T_1)$.
\end{proof1}

\bigskip
Let $\mathcal{R}(T_1)$ denote the set of recolored vertices in the set of labels  $\phi(T_1)$. The above proof states that this set is uniquely defined by the set of edges. Moreover, repeating the proof of Proposition \ref{blaming} implies that for every node $C\ne A$, the blaming vertex $v(\phi(C))$ is also uniquely defined.

\begin{proposition}\label{bdisjoint:property3}
Suppose that $A_0=A,A_1,\ldots,A_{t-1}$ are the nodes of $T_1$. Then in every $\phi(A_i)$ there is a vertex subset $R_i\subset \phi(A_i)$ such that
\begin{enumerate}
  \item $|R_i|\geqslant n-20e\ln n-b$ for any $i=1,\ldots,t-1$;
  \item the sets $R_0,\ldots,R_{t-1}$ have pairwise empty intersection, $R_i\cap R_j=\emptyset$, $i\ne j$;
  \item all the vertices in $R_i$ are colored with dominating color of $A_i$ in the initial coloring $f$;
  \item the vertex $v(\phi(A_i))$ belongs to $R_i$ and is the first vertex of $R_i$ (according to $\sigma$) for every $i>0$.
\end{enumerate}
\end{proposition}
\begin{proof1} Without loss of generality assume that $(A_m,\ldots,A_0)$ is a straight path containing all the bad nodes. Let us define the set $R_i$ for $i=0,\ldots,m$, as follows:
\begin{equation}\label{definition:r_i_1}
    R_i=\{v(\phi(A_i))\} \cup \phi(A_i)\setminus\left(\mathcal{R}(T_1)\cup\bigcup_{j=i+1}^{m} \phi(A_j)\right).
\end{equation}
Here we assume that $\{v(\phi(A_0))\}$ is an empty set. For $i>m$, define
\begin{equation}\label{definition:r_i_2}
    R_i=\{v(\phi(A_i))\} \cup\phi(A_i)\setminus\left(\mathcal{R}(T_1)\cup\bigcup_{F\in T_1\setminus N(A_i)}\phi(F)\right).
\end{equation}
Since every label $\phi(A_i)$ is not degenerate, we have $|\phi(A_i)\cap \mathcal{R}(T_1)|\leqslant 20e\ln n$. For $i>m$, $\phi(A_i)$ is not bad, so $|\phi(A_i)\cap \bigcup_{F\in T_1\setminus N(A_i)}\phi(F)|\leqslant b$. Thus, $|R_i|\geqslant n-20e\ln n-b$. For $i\leqslant m$, the required inequality follows from \eqref{condition:path}.

Suppose that the nodes $A_{m+1},\ldots,A_{t-1}$ are arranged according to the distance from the root in the increasing order. Denote $R'_i=R_i\setminus\{v(\phi(A_i))\}$. For $j<i$,  the set $R'_i$ can have non-empty intersection with $\phi(A_j)$ only if $A_j$ is the parent of $A_i$. But $R'_j$ does not intersect with the labels of children of $A_j$, because all of them belong to $\mathcal{R}(T_1)$. The sets $R'_1,\ldots,R'_m$ do not intersect by the definition \eqref{definition:r_i_1}. If $i\leqslant m<j$ then $A_i$ cannot be a descendant of $A_j$, so $R'_i$ and $R'_j$ do not intersect. Consequently, all the sets $R'_0,\ldots,R'_{t-1}$ do not intersect. But the definitions \eqref{definition:r_i_1} and \eqref{definition:r_i_2} say that all these sets do not intersect with $\mathcal{R}(T_1)$, so addition of non-coinciding vertices $v(\phi(F_i))$ from $\mathcal{R}(T_1)$ does not give any intersecting sets among $R_0,\ldots,R_{t-1}$.

Since $R_i$ do not intersect with the labels of children of $A_i$, all its vertices should be colored with dominating color of $A_i$ in the initial coloring $f$. Moreover, the algorithm says that the vertex $v(\phi(A_i))$ can blame $\phi(A_i)$ only if it is the first non-recolored vertex at the moment when $\phi(A_i)$ is monochromatic. Thus, $v(\phi(A_i))$ is the first vertex in $R_i$.
\end{proof1}

\bigskip
Now we are ready to estimate the probability of the event $\mathcal{B}_2(T_1)$ that $T_1$ is a $b$-disjoint proper h-tree without degenerate labels.
\begin{proposition}\label{bdisjoint:property4}
For any $b$-disjoint proper h-tree $T_1$ of size $t$ without degenerate labels,
\begin{equation}\label{prob:event2}
    {\sf Pr}\left(\mathcal{B}_2(T_1)\right)\leqslant r^{1-n-(n-b)(t-1)}\left(\frac 1{n-20e\ln n-b}\right)^{t-1}(1-p)^{n-20e\ln n}.
\end{equation}
\end{proposition}
\begin{proof1} Due to Proposition \ref{bdisjoint:property2} the probability that for a given dominating color of the root, all the vertices in $\phi(T_1)$ receive appropriate colors is equal to $r^{-m}$, where $m$ is the number of vertices in $\phi(T_1)$. Proposition \ref{bdisjoint:property1} states that $m\geqslant n+(n-b)(t-1)$.

Suppose $A_0$ is the root. Due to Proposition \ref{bdisjoint:property3}, for any $A_i\in T_1$, $A_i\ne A_0$, the vertex $v(\phi(A_i))$ is the first in the set $R_i$ of size at least $n-20e\ln n-b$. All such sets are disjoint, so these events are independent. Hence, the probability of the second ingredient does not exceed $(1/(n-20e\ln n-b))^{t-1}$.

Finally, all the vertices of a special set $R_0\subset \phi(A_0)$ are colored with dominating color $\alpha$ in $f$ and none of them is free. Otherwise, the algorithm would not stop and $\phi(A_0)$ would not be monochromatic in the final coloring. The probability of this ingredient equals $(1-p)^{|R_0|}\leqslant (1-p)^{n-20e\ln n-b}$. It remains to notice that $R_0$ do not intersect with all the sets $R_i$.

All three ingredients of the event $\mathcal{B}_2(T_1)$ are independent. This yields the estimate \eqref{prob:event2}.
\end{proof1}

\bigskip
The last proposition of the paragraph counts the number of h-trees $T$ for which $\phi(T)$ is containing a fixed vertex of $H$. The argument repeats the proof of Proposition 6 from \cite{KozikShab}.

\begin{proposition}\label{bdisjoint:property5}
Let $H=(V, E)$ be a hypergraph with maximum edge degree $\Delta(H)$ and let $v\in\ V$ be its arbitrary vertex. Then the number of h-trees $T$ of size $t$ with $v\in\phi(T)$ is at most $2(4\Delta(H))^{t}$.
\end{proposition}
\begin{proof1}
Let us fix some specific tree structure $S$ of size $t$. We have $t$ possible choices for the node $C$ for which $v\in\phi(C)$ holds and at most $\Delta(H)+1\leqslant 2\Delta(H)$ possible label choices for that node.

Further we extend labelling according to the following rule: for every unlabelled node $Y$ which is adjacent to a labelled node $X$ pick any edge that intersects $\phi(X)$. Each time we have at most $\Delta(H)$ choices for the next node, hence the total number of node labellings constructed in such a way is at most $2\Delta(H)^t$. When all the nodes are labelled, we can uniquely find the vertices that blame the nodes (due to Proposition \ref{bdisjoint:property2}). Hence we obtain an h-tree with given structure $S$.

Clearly every $b$-disjoint h-tree with structure $S$, $v\in\phi(S)$ can be constructed in this way. Therefore $v$ belongs to at most $2t(\Delta(H))^t$ h-trees with structure $S$. The number of possible structures of size $t$ does not exceed $4^t/t$. Hence, the total number of $b$-disjoint h-trees $T$ of size $t$ for which $v\in\phi(T)$ is smaller than $2(4\Delta(H))^t$.
\end{proof1}

\bigskip
Finally, we obtain the following estimate for the local polynomial corresponding to the second bad event.
\begin{align}\label{polynom:event2}
    w^2_v\left(\frac 1{1-\tau_0}\right)&=\sum_{T_1:\;v\in\phi(T_1)}{\sf Pr}(\mathcal{B}_2(T_1))\left(\frac 1{1-\tau_0}\right)^{|{\rm vbl}(\mathcal{B}_2(T_1))|}\leqslant\notag\\
    &\leqslant\sum_{t=1}^{|E|}\sum_{T_1:\;v\in\phi(T_1),\;|T_1|=t}{\sf Pr}(\mathcal{B}_2(T_1))\left(1+\frac 1n\right)^{|{\rm vbl}(\mathcal{B}_2(T_1))|}\leqslant\notag\\
    &\mbox{(using \eqref{prob:event2} and the estimate $|{\rm vbl}(\mathcal{B}_2(T_1))|\leqslant nt$)}\notag\\
    &\leqslant\sum_{t=1}^{|E|}\sum_{T_1:\;v\in \phi(T_1),\;|T_1|=t}r^{1-n-(n-b)(t-1)}\left(\frac 1{n-20e\ln n-b}\right)^{t-1}(1-p)^{n-20e\ln n-b}e^{t}\leqslant\notag\\
    &\mbox{(assuming that $n$ is large and using Proposition \ref{bdisjoint:property5} with the condition \eqref{bound:new})}\notag\\
    &\leqslant\sum_{t=1}^{|E|}2(4\Delta(H))^tr^{1-n-(n-b)(t-1)}\left(\frac 2{n}\right)^{t-1}(1-p)^{n-20e\ln n-b}e^{t}\leqslant\notag\\
    &\leqslant\sum_{t=1}^{|E|}2\left(\frac {4}{(2e)^4}\right)^tn^t\left(\frac 2{n}\right)^{t-1}\left(1-\frac{5\ln n}{n}\right)^{n-20e\ln n}e^{t}\leqslant\notag\\
    &\leqslant n\cdot n^{-5+o(1)}\sum_{t=1}^{|E|}\left(\frac {8e}{(2e)^4}\right)^t=n^{-4+o(1)}\leqslant \frac{1}{10(n+1)}.
\end{align}

\subsubsection{Bad event 3: a large $b$-disjoint proper h-subtree}

Now suppose that the proper h-tree $T_1=O(T)$ is not $b$-disjoint, but there is an h-subtree $T'$ in $T$ such that $T'_1=O(T')$ is $b$-disjoint and has a size at least $\ln n$, $t=|T'_1|\geqslant \ln n$. Let $\mathcal{B}_3(T'_1)$ denote this event. Note that Propositions \ref{bdisjoint:property1}--\ref{bdisjoint:property3} can also be applied to an h-subtree. The main difference with the case of \emph{$b$-disjoint proper h-trees} is the estimate for the probability.
\begin{proposition}\label{subtree:property1}
For any $b$-disjoint proper h-subtree $T'_1$ of size $t$ without degenerate labels,
\begin{equation}\label{prob:event3}
    {\sf Pr}\left(\mathcal{B}_3(T'_1)\right)\leqslant r^{1-n-(n-b)(t-1)}\left(\frac 1{n-20e\ln n-b}\right)^{t-1}.
\end{equation}
\end{proposition}
\begin{proof1}The first two ingredients of the event $\mathcal{B}_2(T_1)$ in the proof of Proposition \ref{bdisjoint:property4} hold for $\mathcal{B}_3(T'_1)$ also. But now we cannot say that a lot of vertices in the $\phi(A_1)$ ($A_1$ is a root of $T_1$) should not be free (in fact, at least one of them should be free), so we omit the third event for $\mathcal{B}_2(T_1)$. The intersection of the first two implies the estimate \eqref{prob:event3}.
\end{proof1}

\bigskip
The number of h-subtrees $T$ of a fixed size for which $\phi(T)$ contains a given vertex can be estimated by using Proposition \ref{bdisjoint:property5}.
Therefore we obtain the following bound for the local polynomial corresponding to the third bad event:
\begin{align}\label{polynom:event3}
    w^3_v\left(\frac 1{1-\tau_0}\right)&=\sum_{T'_1:\;v\in\phi(T'_1)}{\sf Pr}(\mathcal{B}_3(T'_1))\left(\frac 1{1-\tau_0}\right)^{|{\rm vbl}(\mathcal{B}_3(T'_1))|}\leqslant\notag\\
    &\leqslant\sum_{t\geqslant \ln n}\sum_{T'_1:\;v\in\phi(T'_1),\;|T'_1|=t}{\sf Pr}(\mathcal{B}_3(T'_1))\left(1+\frac 1n\right)^{|{\rm vbl}(\mathcal{B}_3(T'_1))|}\leqslant\notag\\
    &\mbox{(using \eqref{prob:event3} and the estimate $|{\rm vbl}(\mathcal{B}_3(T'_1))|\leqslant nt$)}\notag\\
    &\leqslant\sum_{t\geqslant \ln n}\sum_{T'_1:\;v\in\phi(T'_1),\;|T'_1|=t}r^{1-n-(n-b)(t-1)}\left(\frac 1{n-20e\ln n-b}\right)^{t-1}e^{t}\leqslant\notag\\
    &\mbox{(assuming that $n$ is large and using condition \eqref{bound:new})}\notag\\
    &\leqslant\sum_{t\geqslant \ln n}2(4\Delta(H))^tr^{1-n-(n-b)(t-1)}\left(\frac 2{n}\right)^{t-1}e^{t}\leqslant\notag\\
    &\leqslant\sum_{t\geqslant\ln n}2\left(\frac {4}{(2e)^4}\right)^tn^t\left(\frac 2{n}\right)^{t-1}e^{t}\leqslant
    2n\cdot \sum_{t\geqslant \ln n}\left(\frac {8e}{(2e)^4}\right)^t\leqslant\notag\\
    &\leqslant 2n\cdot (2e^3)^{1-\ln n}= n^{-\ln 2-2+o(1)}\leqslant \frac{1}{10(n+1)}.
\end{align}

\subsubsection{Bad event 4: small not $b$-disjoint proper h-subtree}

Let $T$ be an h-tree such that $O(T)$ is not $b$-disjoint. Consider the smallest subtree $Y$ such that $Y'=O(Y)$ is not $b$-disjoint. Let $A$ denote the root of $Y$ and let $F_1,\ldots,F_s$ denote its children in $Y$. Then every subtree $N(F_i)$ satisfies the property that $O(N(F_i))$ is $b$-disjoint, otherwise $Y$ is not the smallest. If the size of $O(N(F_i))$ is greater than $\ln n$ then it is the third bad event and we have already analyzed it. Thus, we may assume that the size of every $O(N(F_i))$ is less than $\ln n$. Since $s\leqslant 20e\ln n$ (recall that we do not have degenerate labels in $T$), the size of $Y'$ is also bounded:
\begin{equation}\label{nobdisjoint:size}
    |Y'|\leqslant \sum_{i=1}^s|O(N(F_i))|+1\leqslant (20e\ln n)\ln n+1\leqslant 30e(\ln n)^2.
\end{equation}
We do not have an equality in \eqref{nobdisjoint:size} since some nodes of $Y$ can be deleted in $O(Y)$ but survive in $O(N(F_i))$ (their copies can belong to $N(F_j)$ for $j\ne i$).

\bigskip
If $Y'$ is not $b$-disjoint then
\begin{itemize}
  \item[(a)] either there is a straight path $C_m,\ldots,C_1,C_0=A$ containing all the bad nodes, but the condition \eqref{condition:path} does not hold;
  \item[(b)] or there is no such a path, so there are two bad nodes $C$ and $D$ such that $C\notin N(D)$ and $D\notin N(C)$.
\end{itemize}
Let $\mathcal{B}_4(Y')$ denote the described event. The following proposition estimate the probability of the event $\mathcal{B}_4(Y')$.

\begin{proposition}\label{badsubtree:property1}
Suppose $Y'$ has a size $t$. Then
\begin{equation}\label{prob:event4}
    {\sf Pr}\left(\mathcal{B}_4(Y')\right)\leqslant r^{1-t(n-bt)}.
\end{equation}
\end{proposition}
\begin{proof1}All the nodes of $Y'$ have different labels, we have deleted all the copies by the operation $O(Y)$. Since $H$ is $b$-simple, every label has at least $n-bt$ vertices which are not contained in any other label. Thus, the total number of vertices is at least $t(n-bt)$. The h-subtree construction guarantees that for a given dominating color of the root, the colors of all the vertices are uniquely defined. The same is true after making the operation $O$ (see Proposition \ref{bdisjoint:property2}). This implies the estimate \eqref{prob:event4}.
\end{proof1}

\bigskip
Now let us estimate the possible number of not $b$-disjoint h-subtrees $T$ for which $\phi(T)$ contains a given vertex of $H$.
\begin{proposition}\label{badsubtree:property2}
The number of not $b$-disjoint h-subtrees $T$ of size $t$ without degenerate labels and containing a given vertex $v\in\phi(T)$ does not exceed
$$
  4\cdot 4^tt^2(\Delta(H))^{t-1}{nt\choose b+1}.
$$
\end{proposition}
\begin{proof1}The tree structure of $T$ can be be chosen in at most $4^t/t$ ways.

Consider the case (a). Recall that here we have a straight path $C_m,\ldots,C_1,C_0$ ($C_0$ is the root) which contains all the bad nodes. But the condition \eqref{condition:path} does not hold. We have to choose three special nodes: a node $D$ with $v\in\phi(D)$, a node $C_m$ and a node $C_j$, $j<m$, for which the condition \eqref{condition:path} fails (in at most $t^3$ ways). Note that $C_m$ is a bad node. A label $\phi(D)$ can be chosen in at most $\Delta(H)+1\leqslant2\Delta(H)$ ways. If $C_m$ does not belong to a straight path connecting $D$ and the root $C_0$, then we can label all the nodes that do not belong to $N(C_m)$ before $C_m$. Just use to the following rule: for every unlabelled node $X$ which is adjacent to a labelled node $Z$ pick any edge that intersects $\phi(Z)$. Each time we have at most $\Delta(H)$ choices for the next label. After that we have to choose $\phi(C_m)$. Since our hypergraph is $b$-simple and $C_m$ is a bad node, it can uniquely defined by some $b+1$ vertices from the set $\bigcup_{F\in Y'\setminus N(C_m)}\phi(F)$. This set is already defined and has a size at most $nt$, so $C_m$ can be labelled in at most ${nt\choose b+1}$ ways. All the remaining nodes in $N(C_m)$ can be labelled by the usual rule, in at most $\Delta(H)$ ways each.

If $C_m$ belongs to a straight path connecting $D$ and the root $C_0$, then we can label all the nodes on the path by usual rule until we reach $C_j$. Due to the complement of \eqref{condition:path} the node $C_j$ can be labelled in at most ${nt\choose b+1}$ ways, because its label $\phi(C)$ should have at least $b+1$ common vertices with already chosen labels $\phi(C_m),\ldots,\phi(C_{j+1})$.

\bigskip
In the case (b) we have two bad nodes $C$ and $F$, which do not lie on the same straight path to the root. Again we have to choose three special nodes: a node $D$ for which $v\in\phi(D)$, a node $C$ and a node $F$. Then either we can label all the nodes that do not belong to $\phi(N(C))$ before labelling $C$ or we can label all the nodes that do not belong to $\phi(N(F))$ before labelling $F$. Indeed, at least one of the nodes $C$ or $F$ does not lie on the straight path from $D$ to the root. Suppose it is $C$. Again we use the usual rule: for every unlabelled node $X$ which is adjacent to a labelled node $Z$ pick any label that intersects $\phi(Z)$. Each time we have at most $\Delta(H)$ choices for the next label. After that we have to choose $C$. Since our hypergraph is $b$-simple and $C$ is a bad node, its label $\phi(C)$ can uniquely defined by some $b+1$ vertices from the set $\bigcup_{F'\in Y'\setminus N(C)}\phi(F')$. This set is already defined and has a size at most $nt$, so $C$ can be labelled in at most ${nt\choose b+1}$ ways. All the remaining nodes in $N(C)$ can be labelled by the usual rule, in at most $\Delta(H)$ ways each.
\end{proof1}

Note that above result gives better bound than Proposition \ref{bdisjoint:property5} if $r$ is big enough in comparison with $n$.

\bigskip
Finally, we obtain the following bound for the local polynomial corresponding to the fourth bad event:
\begin{align}\label{polynom:event4}
    w^4_v\left(\frac 1{1-\tau_0}\right)&=\sum_{Y':\;v\in\phi(Y')}{\sf Pr}(\mathcal{B}_4(Y'))\left(\frac 1{1-\tau_0}\right)^{|{\rm vbl}(\mathcal{B}_4(Y'))|}\leqslant\notag\\
    &\mbox{(using \eqref{nobdisjoint:size})}\notag\\
    &\leqslant\sum_{t\leqslant 30e(\ln n)^2}\sum_{Y':\;v\in\phi(Y'),\;|Y'|=t}{\sf Pr}(\mathcal{B}_4(Y'))\left(1+\frac 1n\right)^{|{\rm vbl}(\mathcal{B}_4(Y'))|}\leqslant\notag\\
    &\mbox{(using \eqref{prob:event4} and the estimate $|{\rm vbl}(\mathcal{B}_4(Y'))|\leqslant nt$)}\notag\\
    &\leqslant\sum_{t\leqslant 30e(\ln n)^2}\sum_{Y':\;v\in\phi(Y'),\;|Y'_1|=t}r^{1-t(n-bt)}e^{t}\leqslant\notag\\
    &\mbox{(assuming that $n$ is large and using condition \eqref{bound:new})}\notag\\
    &\leqslant\sum_{t\leqslant 30e(\ln n)^2}4\cdot4^t(\Delta(H))^{t-1}t^2{nt\choose b+1}r^{1-t(n-bt)}e^{t}\leqslant\notag\\
    &\leqslant\sum_{t\leqslant 30e(\ln n)^2}16\left(\frac {4}{(2e)^4}\right)^{t-1}t^2n^te^{t}(nt)^{b+1}r^{(n-b)(t-1)+1-t(n-bt)}\leqslant
    \notag\\
    &\leqslant\sum_{t\leqslant 30e(\ln n)^2}16e\left(\frac {1}{4e^3}\right)^{t-1}t^2n^t(nt)^{b+1}r^{b(t^2-t+1)+1-n}\leqslant
    \notag\\
    &\mbox{(since $t=O((\ln n)^2)$ and $n$ is large in comparison with $b$)}\notag\\
    &\leqslant\sum_{t\leqslant 30e(\ln n)^2}16e\left(\frac {1}{4e^3}\right)^{t-1} e^{O((\ln n)^3)}r^{O((\ln n)^2)-n}\leqslant
    \notag\\
    &\leqslant 2^{O((\ln n)^3)-n}\leqslant \frac{1}{10(n+1)}.
\end{align}

\subsection{Completion of the proof}
For the application of the Local Lemma \ref{LocalLemma}, we have to check the condition \eqref{local:condition}. For every vertex $v$, we have proved the estimates \eqref{polynom:event1}, \eqref{polynom:event2}, \eqref{polynom:event3}, \eqref{polynom:event4} for the local polynomials. So, their sum $w_v(z)=\sum_{i=1}^4w_v^i(z)$ satisfies the relation
\begin{align*}
   w_v\left(\frac{1}{1-\tau_0}\right)=\sum_{i=1}^4w_v^i\left(\frac{1}{1-\tau_0}\right)\leqslant \frac{4}{10(n+1)}<\frac{1}{n+1}=\tau_0.
\end{align*}
The Local Lemma implies that with positive probability none of the bad events occurs, so the algorithm produces a proper $r$-coloring and the hypergraph $H$ is $r$-colorable. Theorem \ref{thm:main} is proved.

\bigskip
\section{Corollaries}
\subsection{The maximum vertex degree}
The first corollary of the main result is the obvious estimate for the maximum vertex degree in $b$-simple hypergraphs with high chromatic number.
\begin{corollary}\label{corollary:vertexdegree}
If $H$ is a $b$-simple non-$r$-colorable $n$-uniform hypergraph and $n>n_0(b)$ then its maximum vertex degree is at least $1/(2e)^4r^{n-b}$.
\end{corollary}
\begin{proof1}Theorem \ref{thm:main} implies that $H$ contains an edge $A$ with degree at least $1/(2e)^4n\cdot r^{n-b}$. Hence, $A$ contains a vertex of degree at least $1/(2e)^4r^{n-b}$.
\end{proof1}

\subsection{The number of edges}
In \cite{KMRT} Kostochka, Mubayi, R\"odl and Tetali proposed to consider the problem of estimating the minimum possible number of edges in an $n$-uniform $b$-simple hypergraph with chromatic number greater than $r$. Let $m(n,r,b)$ denote the considered extremal value. The authors of \cite{KMRT} showed that for fixed $n$ and $b$, the function $m(n,r,b)$ has the order $\Theta_{n,b}((r\ln r)^{1+1/b})$ as a function of $r$. In this paper we concentrate on the opposite situation: $r$, $b$ are fixed and $n$ grows. In this asymptotic area  Kostochka and Kumbhat \cite{KostKumb} proved that
\begin{equation}\label{bound:kk_bsimple}
  r^{n(1+1/b)}n^{-\varepsilon(n)}\leqslant m(n,r,b)\leqslant c_1\,r^{n(1+1/b)}n^{2(1+1/b)},
\end{equation}
where $\varepsilon(n)>0$ slowly tends to zero as $n\to+\infty$ and $c_1=c_1(b,r)>0$ does not depend on $n$. Later the upper bound in \eqref{bound:kk_bsimple} was improved by Kostochka and R\"odl \cite{KR}, who showed that
\begin{equation}\label{bound:kr_bsimple}
  m(n,r,b)\leqslant c_2\,r^{n(1+1/b)}n^{1+1/b},
\end{equation}
where $c_2=c_2(b,r)>0$ does not depend on $n$. The current best lower bound is due to Kozik \cite{Kozik}:
\begin{equation}\label{bound:kozik_bsimple}
  m(n,r,b)\geqslant \Omega_{r,b}\left(\left(\frac{r^n}{\ln n}\right)^{1+1/b}\right).
\end{equation}

\bigskip
We refine the bound \eqref{bound:kozik_bsimple} as follows.

\begin{corollary}\label{corollary:m(n,r,b)}
For any fixed $r\geqslant 2$, $b\geqslant 2$ and sufficiently large $n>n_0(b)$,
\begin{equation}\label{bound:new_m(n,r,b)}
  m(n,r,b)\geqslant c\cdot r^{n(1+1/b)},
\end{equation}
where $c=c(r,b)>0$ depends only on $r$ and $b$.
\end{corollary}
\begin{proof1}
Here we just follow the argument from \cite{Kozik}. The proof is based on a trimming technique which was proposed by Erd\H{o}s and Lov\'asz and developed by Kostochka and Kumbhat. Suppose that $H=(V,E)$ is an $n$-uniform $b$-simple non-$r$-colorable hypergraph. For every edge $e\in E$, fix a vertex $v(e)$ which has the maximum degree among all the vertices of $e$ (if there are a few such vertices then choose one arbitrarily). Consider the following hypergraph $H'=(V,E')$ where
$$
  E'=\{e\setminus \{v(e)\}:\;e\in E\}.
$$
In other words we remove a vertex with maximum degree from every edge (a trimming procedure).

\bigskip
Let $H_1$ be a hypergraph obtained from $H$ by applying successively the trimming procedure $b$ times. It is clear that $H_1$ is $(n-b)$-uniform, $b$-simple and non-$r$-colorable. So, Corollary \ref{corollary:vertexdegree} implies that $H_1$ contains a vertex $v$ of degree $d$ at least $1/(2e)^4 r^{n-2b}$.

Let $f_1,\ldots,f_d$ denote the edges of $H_1$ containing $v$ and let $e_1,\ldots,e_d$ denote the corresponding edges of $H$. Let $Y$ denote the set of vertices which were removed from the edges containing $v$ in $H$ during the trimming procedure. Since $H$ is $b$-simple every $e_j$ is uniquely defined by some $b$-tuple of $Y$. Therefore
$$
  d\leqslant {|Y|\choose b}\leqslant |Y|^b.
$$
So hypergraph $H$ contains at least $d^{1/b}$ vertices of degree at least $d$. Let $v_1,\ldots,v_m$, $m=\lceil d^{1/b}\rceil$, denote this set of vertices. Also let $d_j$ denote the number of edges which contain $v_j$ and have at most $b-1$ vertices from $v_1,\ldots,v_{j-1}$. Since $H$ is $b$-simple, any other edge containing $v_j$ is uniquely defined by a $b$-tuple from $v_1,\ldots,v_{j-1}$. Thus, $d_j\geqslant d-{j-1\choose b}$ and
$$
  \sum_{j=1}^md_j\geqslant \sum_{j=1}^m\left(d-{j-1\choose b}\right)\geqslant dm-{m\choose b+1}\geqslant d^{1+1/b}-\frac{m^{b+1}}{(b+1)!}\geqslant c_0\,d^{1+1/b},
$$
where $c_0=c_0(b)>0$ depends only on $b$. Recall that we assume that $n$ is sufficiently large in comparison with $b$.

Finally, it remains to note that any edge of $H$ is counted at most $b$ times in the sum $\sum_{j=1}^md_j$. Hence,
$$
  |E|\geqslant \frac 1b \sum_{j=1}^md_j\geqslant \frac{c_0}b\,d^{1+1/b}\geqslant c(r,b)r^{n(1+1/b)}.
$$
\end{proof1}

Note that the lower bound \eqref{bound:new_m(n,r,b)} is only $n^{1+1/b}$ times smaller than the upper bound \eqref{bound:kr_bsimple} for fixed $r$, $b$ and large $n$.

\section{Acknowledgements}

The research of the second author was supported by the program Leading Scientific Schools (grant no. NSh-6760.2018.1). The authors are grateful to Professor Jakub Kozik for the fruitful discussion of the problem.

\renewcommand{\refname}{References}

\end{document}